\documentclass[11pt, twosided]{article}
\usepackage{amsmath}
\usepackage{amsfonts}
\usepackage{color}
\usepackage{amsmath,amsthm,amssymb}
\makeatletter
\newtheorem{thm}{Theorem}[section]
\newtheorem{prop}[thm]{Proposition}
\newtheorem{lem}[thm]{Lemma}

\theoremstyle{definition}
\newtheorem{defi}[thm]{Definition}

\def\proof {{\noindent\bf Proof.}\quad}
\def \endproof{\hfill$\Box$\vspace {3mm}}

\newcommand{\ba}{\begin{array}}
\newcommand{\ea}{\end{array}}
\newcommand{\bt}{\begin{tabular}}
\newcommand{\et}{\end{tabular}}
\newcommand{\btb}{\begin{table}}
\newcommand{\etb}{\end{table}}
\newcommand{\bc}{\begin{center}}
\newcommand{\ec}{\end{center}}
\newcommand{\bea}{\begin{eqnarray}}
\newcommand{\eea}{\end{eqnarray}}
\newcommand{\Bea}{\begin{eqnarray*}}
\newcommand{\Eea}{\end{eqnarray*}}
\newcommand{\beq}{\begin{equation}}
\newcommand{\eeq}{\end{equation}}
\newcommand{\Beq}{\begin{equation*}}
\newcommand{\Eeq}{\end{equation*}}

\newcounter{vcenterstest}
\newlength{\leftlength}
\newlength{\rightlength}
\newlength{\vcentersskip}
\newcommand{\leftcentersright}[4][2]{%
        \settowidth{\leftlength}{#2}%
        \settowidth{\rightlength}{#4}%
        \setcounter{vcenterstest}{#1}%
        \ifthenelse{\value{vcenterstest} = 0}
                {\setlength{\vcentersskip}{0pt}}{}%
        \ifthenelse{\value{vcenterstest} = 1}
                {\setlength{\vcentersskip}{\smallskipamount}}{}%
        \ifthenelse{\value{vcenterstest} = 2}
                {\setlength{\vcentersskip}{\medskipamount}}{}%
        \ifthenelse{\value{vcenterstest} = 3}
                {\setlength{\vcentersskip}{\bigskipamount}}{}%
        \ifthenelse{\value{vcenterstest} = 4}
                {\setlength{\vcentersskip}{1cm}}{}%
        \vskip\vcentersskip
        \noindent#2\hskip-\leftlength%
        \hfill#3\hfill
        \mbox{}\hskip-\rightlength#4%
        \vskip\vcentersskip
        }

\begin{document}
\title{On the Deformation of Lie-Yamaguti algebras}
\author{Jie Lin$^{a}$,\,
  Liangyun Chen$^{b}$\thanks{Corresponding author. {\em E-mail
address:} chenly640@nenu.edu.cn.},\, Yao Ma$^{b}$\\
{\footnotesize\em $^{a}$Sino-European Institute of Aviation
Engineering, Civil Aviation University of China,}\\{\footnotesize\em
 Tianjin 300300, China}\\
{\footnotesize\em $^b$School of Mathematics and Statistics,
Northeast Normal University, Changchun 130024, China}}
\date{}

\maketitle

\baselineskip 18pt
\begin{abstract}
The deformation theory of Lie-Yamaguti algebras is developed by choosing a suitable cohomology. The relationship between the deformation and the obstruction of Lie-Yamaguti algebras is obtained.
\medskip

\end{abstract}

\renewcommand{\theequation}{\arabic{section}.\arabic{equation}}

\section{Introduction}
\qquad In mathematical deformation theory one studies how an object in a certain category of spaces can be varied in dependence of the points of a parameter space. In other words, deformation theory thus deals with the structure of families of objects like varieties, singularities, vector bundles, coherent sheaves, algebras or differentiable maps. Deformation problems appear in various areas of mathematics, in particular in algebra, algebraic and analytic geometry, and mathematical physics.

The mathematical theory of deformations has proved to be a powerful tool in modeling physical reality. For example, (algebras associated with) classical quantum mechanics (and field theory) on a Poisson phase space can be deformed to (algebras associated with) quantum mechanics (and quantum field theory). That is a frontier domain in mathematics and theoretical physics called deformation quantization, with multiple ramifications, avatars and connections in both mathematics and physics. These include representation theory, quantum groups (when considering Hopf algebras instead of associative or Lie algebras), noncommutative geometry and manifolds, algebraic geometry, number theory. The notion of deformation can be applied to a variety of categories that are used to express mathematically the physical reality.

The deformation of algebraic systems has been one of the problems that many mathematical researchers are interested in, since Gerstenhaber studied the deformation theory of algebras in a series of papers [\ref{ref1}]-[\ref{ref5}].  For example, it has been extended to covariant functors from a small category to algebras ([\ref{ref6}]) and to algebraic systems, bialgebras, Hopf algebras ([\ref{ref7}]) by Gerstenhaber and Schack, also to Leibniz pairs and Poisson algebras ([\ref{ref8}]) by Flato, Gerstenhaber and Voronov, to Lie triple systems ([\ref{ref9}]) by F. Kubo and Y. Taniguchi. Inspired by these work, we study the deformation theory of Lie-Yamaguti algebras.

 \begin{defi}$^{[\ref{ref04}]}$ Let $k$ be a field of characteristic zero. A Lie-Yamaguti algebra(LYA for short) is a vector space $T$
over $k$ with a bilinear composition $ab$ and a trilinear composition $[a, b,
c]$ satisfying: \begin{equation} a^{2}=0,\label{eq1}\end{equation}
\begin{equation}[a, a, b]=0, \end{equation} \begin{equation} [a, b, c]+ [b,c,a]+[c,a,b]+(ab)c+(bc)a+(ca)b=0,\end{equation} \begin{equation} [ab, c, d]+[bc, a, d]+[ca, b, d]=0,\end{equation}
\begin{equation} [a, b, cd]=[a, b, c]d+c[a, b, d], \end{equation}
\begin{equation}[a,b,[c,d,e]]=[[a,b,c],d,e]+[c,[a,b,d],e]+[c,d,[a,b,e]],\label{eq6}\end{equation}
for any $a, b, c, d, e \in T$.
\end{defi}
Yamaguti([\ref{ref02}]) called it the general Lie triple system, later it was renamed as Lie triple algebras in [\ref{ref03}]. We will follow here the notation in [\ref{ref04}], and will call this system Lie-Yamaguti algebra.

The LYAs with $xy=0$ for any $x, y$ are exactly the Lie triple systems, which is closely related with symmetric spaces, while the LYAs with $[x, y,  z]=0$ are the Lie algebras. These nonassociative binary-ternary algebras have been treated by several authors in connection with geometric problems on homogeneous spaces [\ref{ref05}]-[\ref{ref09}], and their structure theory has been studied by P. Bentio and A. Elduque in [\ref{ref010}]-[\ref{ref012}]. Less known examples can be found in [\ref{ref010}] where a detailed analysis on the algebraic structure of LYAs arising from homogeneous spaces which are quotients of the compact Lie group $G_2$ is given.

On the other hand, Yamaguti's approach ([\ref{ref013}]) to a cohomology theory of LYA's is intrinsic. His cohomology theory is discussed without going out of a LYA into an enveloping Lie algebra, so that the Yamaguti coboundary is defined in terms of only elements of a LYA.

In this paper, a deformation theory of LYAs will be developed. Specifically, a suitable cohomology group is fund to identify the infinitesimal deformations of LYAs; the analytically rigid LYAs are determined; the obstructions to the integration of an infinitesimal deformation are studied.
 %Due to Gerstenhaber ([\ref{ref2}]), we should pay attention to the following aspects of deformation theory :
%\begin{enumerate}
%\item {A definition of the class of objects within which deformation takes place, and identification of the infinitesimal deformations of a given object with the elements of a suitable cohomology group.}
%\item {A theory of the obstructions to the integration of an infinitesimal deformation.}
%\item {A parameterization of the set of objects. }
%\item {A determination of the natural automorphisms of the parameter space and determination of the rigid objects.}
%\end{enumerate}

\section{Deformation of a Lie-Yamaguti algebra}
\setcounter{equation}{0}
\qquad Let $T$ be a LYA over $k$. Let $k [\![t]\!]=R$ denote the power series ring in
one variable $t$, $K=k((t))$ be the quotient power series field of $R$, and let $T_K$ denote the vector space obtained
from $T$ by extending the coefficient domain from $k$ to $K$, i.e., $T_K=T\otimes_k K.$ Any bilinear function $f: T\times T\rightarrow T$
and trilinear function $g: T\times T\times T\rightarrow T$ (in particular, the multiplications in T) can be extended to a function
 bilinear over $K$ from $T_K\times T_K$ to $T_K$ and a function trilinear from $T_K\times T_K\times T_K$ to $T_K$ respectively. The
  functions $T_K\times T_K$ to $T_K$ and $T_K\times T_K\times T_K$ to $T_K$ which are such two extensions are ``defined over $k$".
  Suppose $f_t: T_K\times T_K\rightarrow T_K$ and $g_t: T_K\times T_K\times T_K\rightarrow T_K$ are two bilinear and trilinear functions
  expressible respectively  in the form \centerline{$f_t(a, b)=F_0(a, b)+tF_1(a, b)+t^2F_2(a, b)+\cdots,$}
  and\centerline{ $g_t(a,b, c)=G_0(a, b, c)+tG_1(a, b, c)+t^2G_2(a, b, c)+\cdots,$}  where every $F_i$ is a bilinear
   function and $G_i$ is a trilinear function defined over $k$ and where we may set $F_0(a, b)=f(a, b)=ab,$ and $G_0(a, b, c)=g(a, b, c)=[a, b, c]$,
    the products in $T$. We denote such a system by $T_t:=(T_K, f_t, g_t).$

According to the deformation theory, $T_t$ is required to be the same kind as $T$, that is, a $k[\![t]\!]$-LYA. Thus the conditions, corresponding
to (\ref{eq1})-(\ref{eq6}),
\beq f_t(a, a)=0,\label{eq7}\eeq
\beq g_t(a, a, b)=0,\label{eq8}\eeq
\beq f_t(f_t(a, b), c)+f_t(f_t(b, c), a)+f_t(f_t(c, a), b)+g_t(a, b, c)+g_t(b, c, a)+g_t(c, a, b)=0,\label{eq9}\eeq
\beq g_t(f_t(a, b), c, d)+g_t(f_t(b, c), a, d)+g_t(f_t(c, a), b, d)=0,\label{eq10}\eeq
\beq g_t(a, b, f_t(c, d))=f_t(g_t(a, b, c), d)+f_t(c, g_t(a, b, d)),\label{eq11}\eeq
\beq g_t(a, b, g_t(c, d, e))=g_t(g_t(a, b, c), d, e)+g_t(c, g_t(a, b, d), e)+g_t(c, d, g_t(a, b, e))\label{eq12}\eeq
must be satisfied. The conditions (\ref{eq7}) and (\ref{eq8}) leads to the obvious equations
\beq F_i(a, a)=0,\label{eq13}\eeq
\beq G_i(a, a, b)=0,\label{eq14}\eeq
for $i\in\mathbb{N}.$ From the conditions (\ref{eq9})-(\ref{eq12}) we have
\beq \sum_{\substack{i+j=n\\ i, j\in\mathbb{N}}}\sum\limits_{(a, b, c)\mbox{cyclic}}F_i(F_j(a, b), c)+\sum\limits_{(a, b, c)\mbox{cyclic}}G_n(a, b, c)=0,\label{eq15}
\eeq
\beq  \sum_{\substack{i+j=n\\ i, j\in\mathbb{N}}}\sum\limits_{(a, b, c)\mbox{cyclic}}G_i(F_j(a, b), c, d)=0,\label{eq16}\eeq
\beq  \sum_{\substack{i+j=n\\ i, j\in\mathbb{N}}}G_i(a, b, F_j(c, d))-F_i(G_j(a, b, c), d)-F_i(c, G_j(a, b, d))=0, \label{eq17}\eeq
\beq   \sum_{\substack{i+j=n\\ i, j\in\mathbb{N}}}G_i(a, b, G_j(c, d, e))
-G_i(G_j(a, b, c), d, e)-G_i(c, G_j(a, b, d), e)-G_i(c, d, G_j(a, b, e))=0\label{eq18}\eeq
for $n\in\mathbb{N}.$ We call these the {\it deformation equations} for a Lie-Yamaguti algebra.

Let $(f_t, g_t), (f'_t, g'_t)$ be deformations of a LYA with $f'_t:=F'_0+tF'_1+t^2F'_2+\cdots $ and\linebreak $g'_t:=G'_0+tG'_1+t^2G'_2+\cdots, $ where $F_0=F'_0=f, G_0=G'_0=g.$ We say that the deformation $(f'_t, g'_t)$ is {\it equivalent} to $(f_t, g_t)$, denoted by $(f'_t, g'_t)\sim (f_t, g_t)$, if there exists a $k[\![t]\!]-$LYA isomorphism $\Phi_t: T'_t\rightarrow T_t$ of the form
$\Phi_t=1_T+t\varphi_1+t^2\varphi_2+\cdots,$ where all $\varphi_i$ are $k-$linear maps $T\rightarrow T$ extended to be $k[\![t]\!]-$linear such that for all $a, b, c\in T,$\\ \centerline{$f'_t(a, b)=\Phi^{-1}_t f_t(\Phi_t(a), \Phi_t(b)):=f_t\ast\Phi_t(a, b),$} \centerline{$g'_t(a, b, c)=\Phi^{-1}_t g_t(\Phi_t(a), \Phi_t(b), \Phi_t(c)):=g_t\ast\Phi_t(a, b, c).$}

When $(F_1, G_1)=(F_2, G_2)=\cdots=0$, we say that $(f_t, g_t)$ is the {\it null deformation} and write $(f_0, g_0).$ A deformation $(f_t, g_t)$ is said to be the {\it trivial deformation} when $(f_t, g_t)\sim (f_0, g_0).$

\section{Cohomology groups of Lie-Yamaguti algebras}
\setcounter{equation}{0} \qquad First we recall the definition of
representation of a LYA [\ref{ref02}].

 Let $\rho$ be a linear mapping of a LYA $T$ into the algebra $E(V)$ of linear endomorphism of a vector space $V$ and $D$ and $\theta$ be the bilinear mappings of $T$ into $E(V)$. $(\rho, D, \theta; V)$ is called a representation of $T$ if $\rho, D$ and $\theta$ satisfy the following relations:
\begin{equation}D(a, b)+\theta(a, b)-\theta(b, a)=[\rho(a), \rho(b)]-\rho(ab),\label{eq19}\end{equation}
 \begin{equation} \theta(a, bc)-\rho(b)\theta(a, c)+\rho(c)\theta(a, b)=0,\label{eq20}\end{equation}
 \begin{equation} \theta(ab, c)-\theta(a, c)\rho(b)+\theta(b, c)\rho(a)=0,\label{eq21}\end{equation}
 \begin{equation} \theta(c, d)\theta(a, b)-\theta(b, d)\theta(a, c)-\theta(a, [b, c, d])+D(b, c)\theta(a, d)=0, \label{eq22}\end{equation}
 \begin{equation} [D(a, b), \rho(c)]=\rho([a, b, c]), \label{eq23}\end{equation}
  \begin{equation} [D(a, b), \theta(c, d)]=\theta([a, b, c], d)+\theta(c, [a, b, d]).\label{eq24}\end{equation}

Following by (\ref{eq19}), we shall sometimes denote by $(\rho, \theta)$ the representation $(\rho, D, \theta)$ simply. From (\ref{eq19}), (\ref{eq20}), (\ref{eq21}), (\ref{eq23}) we have \begin{equation}D(ab, c)+D(bc, a)+D(ca, b)=0.\end{equation}

In a LYA $T$, put $D(a, b): c\mapsto [a, b, c], \theta(a, b): c\mapsto [c, a, b], \rho(a): b\mapsto ab,$ then $(\rho, D, \theta)$ is a representation of $T$ into itself, we call it to be regular. An ideal of $T$ is a subspace of $T$ invariant under this representation.

Let $(\rho, D, \theta; V)$ be a representation of LYA $T$. Let $(f, g)$ be a pair of $2p-$ and $(2p+1)-$linear mappings of $T$ into $V$ such that
 $$f(x_{1}, \cdots, x_{2i-1}, x_{2i}, \cdots, x_{2p})=0$$ and $$g(x_{1}, \cdots, x_{2i-1}, x_{2i}, \cdots, x_{2p+1})=0$$ if $x_{2i-1}=x_{2i}, i=1, 2, \cdots, p.$
 We denote by $C^{n}(T, V), n\geq 1,$ a vector space spanned by such linear mappings. For each element $(f, g)\in C^{2p}(T, V)\times C^{2p+1}(T, V)$, a coboudary operator $\delta: (f, g)\mapsto (\delta_{\textrm{I}}f, \delta_{\textrm{II}}g)$ is a mapping of $C^{2p}(T, V)\times C^{2p+1}(T, V)$ into $C^{2p+2}(T, V)\times C^{2p+3}(T, V)$ defined by the following formulas:
\Bea &&(\delta_{\textrm{I}
}f)(x_{1}, x_{2}, \cdots, x_{2p+2})\\
&=&(-1)^{p}[\rho(x_{2p+1}g(x_{1}, \cdots, x_{2p}, x_{2p+2}))-\rho(x_{2p+2})g(x_{1}, \cdots, x_{2p+1})-g(x_{1}, \cdots, x_{2p}, (x_{2p+1}x_{2p+2}))]\\
&&+\sum\limits_{k=1}^{p}(-1)^{k+1}D(x_{2k-1}, x_{2k})f(x_{1}, \cdots, \hat{x}_{2k-1}, \hat{x}_{2k}, \cdots, x_{2p+2})\\
&&+\sum\limits_{k=1}^{p}\sum\limits_{j=2k+1}^{2p+2}(-1)^{k}f(x_{1}, \cdots, \hat{x}_{2k-1}, \hat{x}_{2k}, \cdots, [x_{2k-1}, x_{2k}, x_{j}], \cdots, x_{2p+2}),\Eea
\Bea &&(\delta_{\textrm{II}
}g)(x_{1}, x_{2}, \cdots, x_{2p+3})\\
&=& (-1)^{p}[\theta(x_{2p+2}, x_{2p+3})g(x_{1},\cdots, x_{2p+1})-\theta(x_{2p+1, x_{2p+3}})g(x_{1}, \cdots, x_{2p}, x_{2p+2})]\\
&&+\sum\limits_{k=1}^{p+1}(-1)^{k+1}D(x_{2k-1}, x_{2k})g(x_{1}, \cdots, \hat{x}_{2k-1}, \hat{x}_{2k}, \cdots, x_{2p+3})\\
&&+\sum\limits_{k=1}^{p+1}\sum\limits_{j=2k+1}^{2p+3}(-1)^{k}g(x_{1}, \cdots, \hat{x}_{2k-1}, \hat{x}_{2k}, \cdots, [x_{2k-1}, x_{2k}, x_{j}], \cdots, x_{2p+3}),\Eea
for $(f, g)\in C^{2p}(T, V)\times C^{2p+1}(T,V), p=1, 2, 3\cdots$\\

In the case $n=1$, we shall only consider a subspace spanned by the diagonal elements $(f, f)\in C^{1}(T, V)\times C^{1}(T, V)$ and $\delta(f, f)=(\delta_{\textrm{I}}f, \delta_{\textrm{II}}f)$ is an element of $C^{2}(T, V)\times C^{3}(T, V)$ defined by $$(\delta_{\textrm{I}}f)(a, b)=\rho(a)f(b)-\rho(b)f(a)-f(ab),$$ $$(\delta_{\textrm{II}}f)(a, b, c)=\theta(b, c)f(a)-\theta(a, c)f(b)+D(a, b)f(c)-f([a, b, c]).$$ Further, for each $(f, g)\in C^{2}(T, V)\times C^{3}(T, V)$ another coboundary operation $\delta^{*}=(\delta^{*}_{\textrm{I}}, \delta^{*}_{\textrm{II}})$ of $C^{2}(T, V)\times C^{3}(T, V)$ into $C^{3}(T, V)\times C^{4}(T, V)$ is defined by \Bea(\delta^{*}_{\textrm{I}}f)(a, b, c)&=&-\rho(a)f(b, c)-\rho(b)f(c, a)-\rho(c)f(a, b)+f(ab, c)+f(bc, a)+f(ca, b)\\
&&+g(a, b, c)+g(b, c, a)+g(c, a, b),\Eea
\Bea (\delta^{*}_{\textrm{II}}g)(a, b, c, d)&=&\theta(a, d)f(b, c)+\theta(b, d)f(c, a)+\theta(c, d)f(a, b)+g(ab, c, d)+g(bc, a, d)\\
&&+g(ca, b, d).\Eea

For each $f\in C^{1}(T, V)$ a direct calculation shows that $\delta_{\textrm{I}}\delta_{\textrm{I}}f=\delta^{*}_{\textrm{I}}\delta_{\textrm{I}}f=0$ and $\delta_{\textrm{II}}\delta_{\textrm{II}}f=\delta^{*}_{\textrm{II}}\delta_{\textrm{II}}f=0.$
Yamaguti showed that for each $(f, g)\in C^{2p}(T, V)\times C^{2p+1}(T, V), \delta_{\textrm{I}}\delta_{\textrm{I}}f=0$ and $\delta_{\textrm{II}}\delta_{\textrm{II}}g=0$, or $\delta\delta(f, g)=0.$ For the case $p\geq 2, Z^{2p}(T, V)\times Z^{2p+1}(T, V)$ is a subspace of $C^{2p}(T, V)\times C^{2p+1}(T, V)$ spanned by $(f, g)$ such that $\delta(f, g)=0$. The cohomology group $H^{2p}(T, V)\times H^{2p+1}(T, V)$ of $T$ associated with representation $(\rho, D, \theta)$ is defined as the factor space $(Z^{2p}(T, V)\times Z^{2p+1}(T, V))/(B^{2p}(T, V)\times B^{2p+1}(T, V)),$ where $B^{2p}(T, V)\times B^{2p+1}(T, V)=\delta(C^{2p-2}(T, V)\times C^{2p-1}(T, V)).$ $H^{1}(T, V)=\{f\in C^{1}(T, V)|\delta_{\textrm{I}}f=0, \delta_{\textrm{II}}f=0\}$ by the above definition. In the case $p=1$, let $Z^{2}(T, V)$ be a subspace of $C^{2}(T, V)$ spanned by $f$ such that $\delta_{\textrm{I}}f=\delta^{*}_{\textrm{I}}f=0$ and $Z^{3}(T, V)$ be a subspace of $C^{3}(T, V)$ spanned by $g$ such that $\delta_{\textrm{II}}g=\delta^{*}_{\textrm{II}}g=0,$ then $H^{2}(T, V)\times H^{3}(T, V)$ is defined as the factor space $(Z^{2}(T, V)\times Z^{3}(T, V))/(B^{2}(T, V)\times B^{3}(T, V))$, where $B^{2}(T, V)\times B^{3}(T, V)=\{\delta(f, f)|f\in C^{1}(T, V)\}.$ Let $f$ be a linear mapping of $T$ into a representation space $V$, $f$ is called a derivation of $T$ into $V$ if $f(ab)=\rho(a)f(b)-\rho(b)f(a)$ and $f([a, b, c])=\theta(b, c)f(a)-\theta(a, c)f(b)+D(a, b)f(c).$ If $(\rho, \theta)$ is the regular representation of $T$, then $f$ is a derivation of $T$. Then, by the definition of $H^{1}(T, V)$, $H^{1}(T, V)$ is the vector space spanned by derivations of $T$ into $V$.

\section{Infinitesimal of deformation}
\setcounter{equation}{0}
\qquad Let us return to the deformation equations (\ref{eq13})-(\ref{eq18}). It follows from (\ref{eq13}) and (\ref{eq14}) that each $(F_i, G_i)$ can be viewed as an elements of the product space $C^2(T, T)\times C^3(T, T)$. If we take $n=1$ in (\ref{eq15})-(\ref{eq18}), then we have $\delta^*_\textrm{I} F_1=0, \delta^*_{\textrm{II}}G_1=0, \delta_\textrm{I}F_1=0$ and $\delta_{\textrm{II}}G_1=0$ respectively. For $n\geq 2,$ (\ref{eq15})-(\ref{eq18}) can be expressed as:
\beq \sum\limits_{i=1}^{n-1}[F_i(F_{n-i}(a, b), c)+F_i(F_{n-i}(b, c), a)+F_i(F_{n-i}(c, a), b)]=-\delta^*_{\textrm{I}}F_n(a, b, c),\label{eq25}  \eeq
\beq \sum\limits_{i=1}^{n-1}[G_i(F_{n-i}(a, b), c, d)+G_i(F_{n-i}(b, c), a, d)+G_i(F_{n-i}(c, a), b, d)]=-\delta^*_{\textrm{II}}G_n(a, b, c, d),\label{eq26}\eeq
\beq \sum\limits_{i=1}^{n-1}[G_i(a, b, F_{n-i}(c, d))-F_i(G_{n-i}(a, b, c), d)-F_i(c, G_{n-i}(a, b, d))]=-\delta_{\textrm{I}}F_n(a, b, c, d),\label{eq27}\eeq
\bea &&\sum\limits_{i=1}^{n-1}[G_i(a, b, G_{n-i}(c, d, e))-G_i(G_{n-i}(a, b, c), d, e)-G_i(c, G_{n-i}(a, b, d), e)\\ \nonumber &&-G_i(c, d, G_{n-i}(a, b, e))]
 =-\delta_{\textrm{II}}G_n(a, b, c, d, e).\label{eq28}\eea
The infinitesimal of the deformation $(f_t, g_t)$ is $(F_1, G_1)$. Since $(F_1, G_1)\in Z^2(T, T)\times Z^3(T, T)$, this concept suits the aspects of the Gerstenhaber's deformation theory.

Assume that deformations $(f_t, g_t)$ and $(f'_t, g'_t)$ are equivalent under $\Phi_t=1_T+t\varphi_1+t^2\varphi_2+\cdots,$ where \centerline{$f_t:=F_0+tF_1+t^2F_2+\cdots, f'_t:=F'_0+tF'_1+t^2F'_2+\cdots,$} and \centerline {$g_t:=G_0+tG_1+t^2G_2+\cdots, g'_t:=G'_0+tG'_1+t^2G'_2+\cdots$ }
with $F_0=F'_0=f, G_0=G'_0=g.$ The defining equations $f'_t=f_t\ast\Phi_t$, $g'_t=g_t\ast\Phi_t,$ i.e., \\ \centerline{$\Phi_t(f'_t(a, b))= f_t(\Phi_t(a), \Phi_t(b)),$ $\Phi^{-1}_t(g'_t(a, b, c))= g_t(\Phi_t(a), \Phi_t(b), \Phi_t(c))$} are equivalent to \quad $\sum\limits_{i+j=n}\varphi_j\circ F'_i=\sum\limits_{i+j=n}F_i\circ \tilde{\Phi}_j,\quad \sum\limits_{i+j=n}\varphi_j\circ G'_i=\sum\limits_{i+j=n}G_i\circ \tilde{\Psi}_j,$\\ or
 \beq F'_n=F_n+\sum\limits_{i=0}^n (F_i\circ \tilde{\Phi}_{n-i}-\varphi_{n-i}\circ F'_i),\quad G'_n=G_n+\sum\limits_{i=0}^n (G_i\circ \tilde{\Psi}_{n-i}-\varphi_{n-i}\circ G'_i),\label{eq29}\eeq where $\varphi_0=1_T$ and \beq F_i\circ \tilde{\Phi}_j(a, b)=\sum\limits_{k+l=j}F_i(\varphi_k(a),\, \varphi_l(b)), G_i\circ \tilde{\Psi}_j(a, b, c)=\sum\limits_{k+l+m=j}G_i(\varphi_k(a), \varphi_l(b), \varphi_m(c)).\label{eq30}\eeq
 For $n=1,$ one has $F'_1-F_1=\delta_{\textrm{I}}\varphi_1 $ and $G'_1-G_1=\delta_{\textrm{II}}\varphi_1.$

 Thus we have the following theorem:
 \begin{thm}
 {Let $(f_t, g_t), (f'_t, g'_t)$ be equivalent deformations of a Lie-Yamaguti algebra $(T, f, g)$, then the first-order terms of them belong to the same cohomology class in the cohomology group $H^2(T, T)\times H^3(T, T).$}
 \end{thm}
\section{Rigidity}

\qquad A LYA $T$ is {\it analytically rigid} if every deformation $(f_t, g_t)$ is equivalent to the null deformation $(f_0, g_0)$. As the deformation theory of algebras with binary products, such as associative algebras and Lie algebras, and that with ternary products, such as Lie triple systems, we have a fundamental theorem.
\begin{thm}
{If $T$ is a Lie-Yamaguti algebra with $H^2(T, T)\times H^3(T, T)=0,$ then $T$ is analytically rigid.}
\end{thm}
\proof Let $(f_t, g_t)$ be a deformation of a LYA $(T, f, g)$ with $f_t=f+t^rF_r+t^{r+1}F_{r+1}+\cdots,$ $g_t=g+t^rG_r+t^{r+1}G_{r+1}+\cdots,$ i.e., $F_1=F_2=\cdots=F_{r-1}=0, G_1=G_2=\cdots=0.$ It follows from (\ref{eq25})-(\ref{eq28})
that $\delta(F_r, G_r)=\delta^*(F_r, G_r)=0,$ i.e., $(F_r, G_r)\in Z^2(T, T)\times Z^3(T, T).$ By our assumption $H^2(T, T)\times H^3(t, T)=0$, we can find $\alpha_r\in C^1(t, T)$ such that $(F_r, G_r)=\delta(\alpha_r, \alpha_r).$  Now consider the deformation $(f'_t, g'_t)$ with $f'_t=f_t\ast(1_T-t^r\alpha_r), g'_t=g_t\ast(1_T-t^r\alpha_t).$ In this case, Eq. (\ref{eq29}) is\\
 \centerline {$F'_r=F_r+f\circ\tilde{\Phi}_r-(-\alpha_r)\circ f=F_r-\delta_{\textrm{I}}\alpha_r=0, G_r+g\circ\tilde{\Psi}_r-(-\alpha_r)\circ g=G_r-\delta_\textrm{II}\alpha_r=0.$}
 Hence \centerline{$f'_t=f+t^{r+1}F'_{r+1}+\cdots,\quad  g'_t=g+t^{r+1}G'_{r+1}+\cdots.$} By induction, one can prove $(f_t, g_t)\sim (f_0, g_0).$
 \endproof
 \section{Obstruction, integration}
 \begin{defi}
 A cocycle $(f_1, g_1)\in Z^2(T, T)\times Z^3(T, T)$ is said to be {\it integrable} if there exists a pair of one parameter family $(f_t, g_t)$ such that $(f_1, g_1)$ is the first-order term, in other words, $f_t=f+t^1F_1+t^2F_2+\cdots, g_t=g+t^1G_1+t^2G_2+\cdots.$
 \end{defi}
  Now let us return to the deformation equations (\ref{eq27}) and (\ref{eq28}). Suppose that we have already had $(F_1, G_1), \cdots, (F_{n-1}, G_{n-1}).$
   We want to find $(F_n, G_n)$ satisfying (\ref{eq27}) and (\ref{eq28}). But there is an obstruction to do so. The verification of this fact sets us a
    long computation, and as a result convinces us that the cohomology we choose here is a suitable one for the deformation theory of LYA. For $i\in[\![1,n-1]\!],$ we denote\\
   \centerline{$F_i\star G_{n-i}: (a, b, c, d)\mapsto G_i(a, b, F_{n-i}(c, d))-F_i(G_{n-i}(a, b, c), d)-F_i(c, G_{n-i}(a, b, d)),$}
  \centerline{ $ \begin{array}{cccl}
      G_i\bigtriangleup G_{n-i}: & (a, b, c, d, e) & \mapsto &  G_i(a, b, G_{n-i}(c, d, e))-G_i(G_{n-i}(a, b, c), d, e) \\
       &  &  & -G_i(c, G_{n-i}(a, b, d), e)-G_i(c, d, G_{n-i}(a, b, e)).
    \end{array}$}
  \begin{lem}
  Suppose $ i\in[\![1, n-1]\!].$ If $(F_i, G_i)\in Z^2(T, T)\times Z^3(T, T),$ then\\
   \centerline {$(F_i\star G_{n-i}+F_{n-i}\star G_i, G_i\bigtriangleup G_{n-i}+G_{n-i}\bigtriangleup G_i)\in Z^4(T, T)\times Z^5(T, T).$}
  \end{lem}
  \proof
  The lemma follows from the following formulas:\\
  $\begin{array}{ll}
     &\delta_{\textrm{I}}(F_i\star G_{n-i})(x_1, x_2, x_3, x_4, x_5, x_6)  \\
     = & -\delta_{\textrm{I}}F_i(x_1, x_2, x_5, G_{n-i}(x_3, x_4, x_6))+\delta_{\textrm{I}}F_i(G_{n-i}(x_1, x_2, x_3), x_4, x_5, x_6) \\
      & +\delta_{\textrm{I}}F_i(x_3, G_{n-i}(x_1, x_2, x_4), x_5, x_6)+\delta_{\textrm{I}}F_i(x_3, x_4, x_5, G_{n-i}(x_1, x_2, x_6)) \\
      & -\delta_{\textrm{I}}F_i(x_1, x_2, G_{n-i}(x_3, x_4, x_5), x_6)+\delta_{\textrm{I}}F_i(x_3, x_4, G_{n-i}(x_1, x_2, x_5), x_6) \\
      & +\delta_{\textrm{II}}G_i(x_1, x_2, x_3, x_4, F_{n-i}(x_5, x_6))-G_i(x_1, x_2, \delta_{\textrm{I}}F_{n-i}(x_3, x_4, x_5, x_6)) \\
      & -F_i(x_5, \delta_{\textrm{II}}G_{n-i}(x_1, x_2, x_3, x_4, x_6))-F_i(\delta_{\textrm{II}}G_{n-i}(x_1, x_2, x_3, x_4, x_5), x_6) \\
      & +G_i(x_3, x_4, \delta_{\textrm{I}}F_{n-i}(x_1, x_2, x_5, x_6))+G_{n-i}(x_3, x_4, x_6)G_i(x_1, x_2, x_5) \\
      & -G_i(x_3, x_4, x_6)G_{n-i}(x_1, x_2, x_5)+G_i(x_3, x_4, x_5)G_{n-i}(x_1, x_2, x_6)\\
      &-G_{n-i}(x_3, x_4, x_5)G_i(x_1, x_2, x_6)-[G_{n-i}(x_1, x_2, x_3), x_4, F_i(x_5, x_6)] \\
      & +[G_i(x_1, x_2, x_3), x_4, F_{n-i}(x_5, x_6)]+[G_{n-i}(x_1, x_2, x_4), x_3, F_i(x_5, x_6)]\\
      &-[G_i(x_1, x_2, x_4), x_3, F_{n-i}(x_5, x_6)],
   \end{array}
  $\\
  $\begin{array}{ll}
     &\delta_{\textrm{II}}(G_i\bigtriangleup G_{n-i})(x_1, x_2, x_3, x_4, x_5, x_6, x_7)\\
    =&-\delta_{\textrm{II}}G_i(x_1, x_2, G_{n-i}(x_3, x_4, x_5), x_6, x_7)+\delta_{\textrm{II}}G_i(G_{n-i}(x_1, x_2, x_3), x_4, x_5, x_6, x_7)\\
     &+\delta_{\textrm{II}}G_i(x_3, G_{n-i}(x_1, x_2, x_4), x_5, x_6, x_7)+\delta_{\textrm{II}}G_i(x_1, x_2, G_{n-i}(x_3, x_4, x_6), x_5, x_7) \\
     &+ \delta_{\textrm{II}}G_i(x_3, x_4, G_{n-i}(x_1, x_2, x_5), x_6, x_7)-\delta_{\textrm{II}}G_i(x_3, x_4, G_{n-i}(x_1, x_2, x_6), x_5, x_7)\\
     &+\delta_{\textrm{II}}G_i(x_1, x_2, x_3, x_4, G_{n-i}(x_5, x_6, x_7))-\delta_{\textrm{II}}G_i(x_1, x_2, x_5, x_6, G_{n-i}(x_3, x_4, x_7)) \\
     &+\delta_{\textrm{II}}G_i(x_3, x_4, x_5, x_6, G_{n-i}(x_1, x_3, x_7))-G_i(\delta_{\textrm{II}}G_{n-i}(x_1, x_2, x_3, x_4, x_5), x_6, x_7) \\
     &-G_i(x_5, \delta_{\textrm{II}}G_{n-i}([x_1, x_2, x_3], x_4, x_6), x_7)-G_i(x_5, x_6, \delta_{\textrm{II}}G_{n-i}(x_1, x_2, x_4, x_7)) \\
     &+G_i(x_3, x_4, \delta_{\textrm{II}}G_{n-i}(x_1, x_2, x_5, x_6, x_7))-G_i(x_1, x_2, \delta_{\textrm{II}}G_{n-i}(x_3, x_4, x_5, x_6, x_7)) \\
     &+[G_i(x_1, x_2, x_6), G_{n-i}(x_3, x_4, x_5), x_7]-[G_{n-i}(x_1, x_2, x_6), G_i(x_3, x_4, x_5), x_7] \\
     &-[G_{n-i}(x_1, x_2, x_3), x_4, G_i(x_5, x_6, x_7)]+[G_i(x_1, x_2, x_3), x_4, G_{n-i}(x_5, x_6, x_7)]\\
     &+[G_{n-i}(x_3, x_4, x_6), x_5, G_i(x_1, x_2, x_7)]-[G_i(x_3, x_4, x_6), x_5, G_{n-i}(x_1, x_2, x_7)]\\
     &-[G_i(x_3, x_4, x_6), G_{n-i}(x_1, x_2, x_5), x_7]+[G_{n-i}(x_3, x_4, x_6), G_i(x_1, x_2, x_5), x_7]\\
     &-[G_i(x_1, x_2, x_5), x_6, G_{n-i}(x_3, x_4, x_7)]+[G_{n-i}(x_1, x_2, x_5), x_6, G_i(x_3, x_4, x_7)]\\
     &+[G_i(x_3, x_4, x_5), x_6, G_{n-i}(x_1, x_2, x_7)]-[G_{n-i}(x_3, x_4, x_5), x_6, G_i(x_1, x_2, x_7)]\\
     &+[G_{n-i}(x_1, x_2, x_4), x_3, G_i(x_5, x_6, x_7)]-[G_i(x_1, x_2, x_4), x_3, G_{n-i}(x_5, x_6, x_7)]\\
     &-[G_{n-i}(x_1, x_2, x_6), x_5, G_i(x_3, x_4, x_7)]+[G_i(x_1, x_2, x_6), x_5, G_{n-i}(x_3, x_4, x_7)].\qquad\Box
   \end{array}$\\
The following result is a consequence of this lemma.
\begin{prop}
Let $(f_t, g_t)$ be a deformation of a Lie-Yamaguti algebra $(T, f, g),$ where\linebreak $f_t=f+t^1F_1+t^2F_2+\cdots, g_t=g+t^1G_1+t^2G_2+\cdots.$ Then\\
 \centerline{$(\sum\limits_{i=1}^{n-1}F_i\star G_{n-i}, \sum\limits_{i=1}^{n-1}G_i\bigtriangleup G_{n-i})\in Z^4(T, T)\times Z^5(T,T).$}
\end{prop}
 The following is the third fundamental theorem.
 \begin{thm}
 If $T$ is a Lie-Yamaguti algebra with $H^4(T, T)\times H^5(T, T)=0,$ then every element in $Z^2(T, T)\times Z^3(T, T)$ is integrable.
 \end{thm}
\section{Acknowledgements}
 The first author gratefully acknowledges the support of  NSFC (No. 11226054) and Scientific Research Foundation of Civil Aviation University of China (No. 09QD08X).
  The second author gratefully acknowledges the support of   NNSF of China (No. 11171057),  Natural
Science Foundation of  Jilin province (No. 201115006) and Scientific
Research Foundation for Returned Scholars
    Ministry of Education of China.

\end{document}